\documentclass[12pt,leqno]{article}

\usepackage{amssymb}
\usepackage{amsmath}

\def\oo{{\omega}}
\def\ooo{{\Omega}}

\def\vp{{\varepsilon}}

\def\pp{{\partial}}
\def\dd{\displaystyle}
\def\vf{{\varphi}}
\def\lbb{{\lambda}}

\def\calo{\mathcal{O}}
\def\calf{\mathcal{F}}

\def\barr{\begin{array}}
\def\earr{\end{array}}
\def\ov{\overline}
\def\wt{\widetilde}
\def\vspp{\vspace*{1,5mm}\\ }

\def\ct{continuous}
\def\ff{\forall }
\def\fwg{following}
\def\3{\subset }
\def\1{^{-1}}

\def\rr{\mathbf{R}}

\def\9{{\infty}}
\def\<{\left<}
\def\>{\right>}

\def\({\left(}
\def\){\right)}

\def\D{\Delta }

\def\bld#1#2{{\buildrel{#1}\over{#2}}}
\def\st#1#2{{\mathrel{\mathop{#2}\limits_{#1}}{}\!}}
\def\q{\quad }
\def\na{{\nabla}}
\def\wt{\widetilde}
\def\de{{\delta}}
\def\vf{{\varphi}}
\def\oo{{\omega}}
\def\ooo{{\Omega}}
\def\pp{{\partial}}
\def\barr{\begin{array}}
\def\earr{\end{array}}
\def\calo{{\cal O}}
\def\dd{\displaystyle}

\def\sk{\smallskip}

\def\pas{\mathbb{P}\mbox{-a.s.}}
\def\b{{\beta}}
\def\eq{equa\-tion}
\def\fc{function}
\def\dd{\displaystyle}
\def\vsp{\vspace*{2mm}\\ }

\def\({\left(}
\def\){\right)}
\def\<{\left<}
\def\>{\right>}
\def\rr{{\mathbb{R}}}
\def\nn{{\mathbb{N}}}
\def\dubluunu{\mbox{$1\!\!\!\,\rule{0,1mm}{2,4mm}\,$}}

\newtheorem{theorem}{Theorem}[section]
\newtheorem{lemma}[theorem]{Lemma}

\newtheorem{definition}[theorem]{Definition}
\newtheorem{remark}[theorem]{Remark}

\usepackage{amssymb,amsfonts,amsmath}
\usepackage{amsfonts}
\usepackage[usenames]{color}

\title{\bf Stochastic porous media equations and
self-organized criticality: convergence to the critical state\\  in all dimensions}

\author{{\bf Viorel Barbu}$^1$\thanks{Supported by the DFG through SFB 701 and IRTG 1132 as well as by the BiBoS Research Centre.}, {\bf Michael R\"ockner}$^2$\\
$^1$ {\small Octav Mayer Institute of Mathematics of Romanian Academy,}\\ {\small 700500 Ia\c si, Romania}\\ $^2$ {\small Fakultat f\"ur Mathematik,
Universitat Bielefeld,}\\ {\small  Postfach 100131, 33501
Bielefeld, Germany}}
\date{}

\begin{document}

\maketitle

\begin{abstract} \noindent If $X=X(t,\xi)$ is the
solution to the stochastic porous media equation in
$\calo\subset\rr^d$, $1\le d\le 3,$ modelling the self-organized
cri\-ti\-ca\-lity \cite{5} and $X_c$ is the critical state, then
it is proved that  $\int^\9_0m(\calo\setminus\calo^t_0)dt<\9,$
$\pas$
 and
 $\lim_{t\to\9}\int_\calo|X(t)-X_c|d\xi=\ell<\9,\ \pas$
 Here, $m$ is the Lebesgue measure and $\calo^t_c$ is the critical
 region
 $\{\xi\in\calo;$ $ X(t,\xi)=X_c(\xi)\}$ and $X_c(\xi)\le X(0,\xi)$ a.e. $\xi\in\calo$. If the stochastic Gaussian perturbation has only finitely many modes (but is still function-valued), $\lim_{t\to\9}\int_K|X(t)-X_c|d\xi=0$ exponentially fast for all compact $K\subset\calo$ with probability one, if the noise is sufficiently strong. We also recover that in the deterministic case   $\ell=0$.\sk\\
{\bf Key words and phrases:} porous media equation, multiplicative noise, self-organized criticality, Ito formula.\sk\\
{\bf2000 Mathematics Subject Classification:} 76S05; 60H15.\end{abstract}

\section{Introduction}

The self-organized criticality is a property of dynamical  systems
which have a critical point as an attractor and which emerges
spontaneously to this attractor. If $X=X(t,\xi)$, $t\ge0,$
$\xi\in\calo\subset\rr^d$, $d=1,2,3,$ is the state of the system
distributed in the spatial domain $\calo$ and if $X_c=X_c(\xi)$ is
a critical state, then $X(t,\cdot)$ divides the space into the
\fwg\ three spatial regions:\medskip

\begin{tabular}{rclcl}
{\it critical region}&&$\calo^t_c$&=&$\{\xi\in\calo;\
X(t,\xi)=X_c(\xi)\},$\vspp {\it subcritical
region}&&$\calo^t_-$&=&$\{\xi\in\calo;\ X(t,\xi)<X_c(\xi)\},$\vspp
{\it supercritical region}&&$\calo^t_+$ &=&$\{\xi\in\calo;\
X(t,\xi)>X_c(\xi)\},$\end{tabular}\medskip

\noindent The main feature of the self-criticality  phenomena is
that the subcritical and supercritical regions are unstable and
absorbed in time by the critical region via an autonomous
mechanism. The standard model of self-organized criticality is the
celebrated {\it sand-pile} model introduced by Bak, Tang and
Wiesenfeld \cite{1}, which is formalized via automation theory
(\cite{2}) and leads to parabolic nonlinear equations of porous
media type
\begin{equation}\label{e1.1}
\frac{\pp X}{\pp t}=a\D H(X-X_c)\mbox{\ \ in\ \
}(0,\9)\times\calo,
\end{equation}
where $a>0$ and $H$ is the Heaviside functions.  (See, also,
\cite{3} for a complete description of this model.) In the
presence of a stochastic Gaussian perturbation, the model is best
described by the stochastic (porous media) equation
\begin{equation}\label{e1.2}
\barr{l}
dX(t)-a\D H(X(t)-X_c)=\sigma(X(t)-X_c)dW_t\ \ \mbox{in }(0,\9)\times\calo,\vspp
X(0)=x\ \ \mbox{in }\calo.\earr
\end{equation}
In \cite{5}, existence and uniqueness of solutions to \eqref{e1.2}
are shown and it is also proved that we have finite-time extinction
of $t\to X(t)-X_c$ with positive probability in $1-D$. In terms of
 self-organized criticality behavior, this means that the
subcritical and supercritical regions are absorbed in finite-time
with positive probability by the critical region $\calo^t_c$. Our
aim here is to establish a similar result in dimensions $d=2,3$, at least asymptotically.
The first main result, Theorem \ref{t2.1} below, amounts to saying that
"for almost all $\{t_n\}\to\9$" we have
\begin{equation}\label{e1.3}
\lim_{n\to\9}m(\calo\setminus\calo^{t_n}_c)=0,\ \pas,
\end{equation}
where $m$ is the Lebesgue measure on $\calo$ and we assume for the
initial state that $x\ge X_c$ a.e. in $\calo$. The second main result, Theorem \ref{t2.2} below, says that $X(t)$, multiplied by the exponential of the function-valued noise, converges to $X_c$ in $L^1(\calo)$ asymptotically and that, if the noise is nondegenerate away from the boundary of $\calo$ (see \eqref{e2.10a} below), then $X(t)$ itself converges asymptotically to $X_c$ locally in $L^1(\calo)$ exponentially fast.

\bigskip\noindent{\bf Notation.} In the \fwg, $\calo$ is a bounded
and open subset of $\rr^d$, $d\ge1$, with smooth boundary, and $L^p(\calo)$, $1\le
p\le\9$, is the space of all $p$-integrable functions in$\calo$
with the usual norm denoted by $|\cdot|_p$. For $k=1,2$,
$H^k(\calo)$, $H^1_0(\calo)$ and $H\1(\calo)$ are standard Sobolev
spaces on $\calo$. More precisely, $H^1_0(\calo)$ is the subspace
of functions $u\in H^1(\calo)$ with zero trace on the boundary
$\pp\calo$ of $\calo$ and $H=H\1(\calo)$ is the dual of
$H^1_0(\calo)$ with the norm
$$|u|_{-1}=\left<A\1u,u\right>_2^{\frac12}.$$
Here, $A=-\D,\ D(A)=H^1_0(\calo)\cap H^2(\calo)$ and
$\left<\cdot,\cdot\right>_2$ is the scalar product of~$L^2(\calo)$.

Everywhere in the \fwg, $\{\ooo,\calf,\calf_t,\mathbb{P}\}$ is a
stochastic basis and $\{\b_j\}^\9_{j=1}$ is a sequence of mutually
independent Brownian motions which induces the filtration
$\{\calf_t\}_{t\ge0}$. By $L^q(0,T;L^p(\ooo,Y))$, where $Y$ is a
Hilbert space, we denote the space of all $q$-integrable processes
$u:(0,T)\to L^p(\ooo,Y)$. By $C([0,T];L^2(\ooo,Y))$ we denote the
space of all $Y$-valued   processes which are mean-square
continuous on $[0,T]$.

\section{Hypotheses and the main result}
\setcounter{equation}{0}

Consider the equation
\begin{equation}\label{e2.1}
\barr{ll} dX(t)-a\D\psi(X(t))dt\ni\sigma(X(t))dW_t\mbox{\ \ in
}(0,\9)\times\calo,\vsp X(0,\xi)=x(\xi),\ \ \xi\in\calo,\vsp
\psi(X(t,\xi))\ni0,\ \ \mbox{on }(0,\9)\times\pp\calo. \earr
\end{equation}
Here, $a$ is a positive constant and
\bigskip

(H1) $\psi(r)={\rm sign}\ r,$\bigskip

\noindent where ${\rm sign}\ r=r|r|^{-1}$ if $r\ne0$, ${\rm sign}\
0=[-1,1]$,\bigskip

(H2) $\sigma(X)dW=\dd\sum^\9_{k=1}\mu_k X e_k d\b_k, $\bigskip

\noindent where  $\{\mu_k\}$ is a sequence of real numbers such
that
\begin{equation}\label{e2.2}
\sum^\9_{k=1}\mu^2_k\lbb^2_k<\9\end{equation} and $\{e_k\}$ is the
orthonormal basis in $L^2(\calo)$ consisting of eigenvectors of
$A$ with eigenvalues $\{\lbb_k\}$, that is, $Ae_k=\lbb_k e_k,$
$k=1,\!...$.  Here $\{\lbb_j\}^\9_{j=1}$ is taken in
increasing order.

\begin{definition}\label{d2.1} {\rm Let $x\in H=H\1(\calo)$. An
$H$-valued continuous $\calf_t$ adapted process $X=X(t)$ is said
to be a solution to \eq\ \eqref{e2.1} if, on every interval $(0,T),\ T>0$,
\begin{equation}\label{e2.3}
X\in L^1(\ooo\times(0,T)\times\calo)\cap
L^2(0,T;L^2(\calo,H))\end{equation} and there is $\eta\in
L^1(\ooo\times(0,T)\times\calo)$ such that
\begin{equation}\label{e2.4}
\barr{ll}
\<X(t),e_j\>_2=\!\!\!&\<x,e_j\>_2+\dd\int^t_0\int_\calo\eta(s,\xi)\D
e_j(\xi)d\xi\,ds\vsp
&+\dd\sum^\9_{k=1}\mu_k\int^t_0\<X(s)e_k,e_j\>_2 d\b_k(s),\ \ff
j\in\nn,\ t\in[0,T],\earr\end{equation}
\begin{equation}\label{e2.5}
\eta\in\psi(X)\mbox{\ \ a.e. on
}(0,T)\times\calo\times\ooo.\end{equation} }\end{definition}

One of the main results established in  \cite{5} (see, also,
\cite{4}) is that, for each $x\in L^p(\calo)$, $p\ge4$, {\it there
is a unique solution $X\in L^\9(0,T;L^p(\ooo,\calo))\cap
L^2(\ooo,C([0,T];H\1(\calo)))$ to equation \eqref{e2.1}. Moreover,
if $x\ge0$  a.e. in $\calo$, then $X\ge0$ a.e., $\pas$} (See
\cite[Theorem 2.2]{5}.) Other existence results for the stochastic porous media equation (\ref{e2.1}) for general maximal monotone functions $\psi$ with  the range all of $R$   were established in \cite{6}, \cite{8}.

In this paper we prove the \fwg\ asymptotic  results for solutions
to \eq~\eqref{e2.1}.

\begin{theorem}\label{t2.1} Assume that {\rm(H1)} and {\rm(H2)}
hold and that $x\in L^4(\calo)$, $x\ge0$, on $\calo.$ Then
\begin{equation}\label{e2.6}
\lim_{t\to\9}\int_\calo X(t,\xi)d\xi=\ell<\9,\ \pas\end{equation}
and
\begin{equation}\label{e2.7}
\int^\9_0m(\calo\setminus\calo^t_0)dt<\9,\ \pas,\end{equation}
where $m$ is the Lebesgue measure and
$$\calo^t_0=\{\xi\in\calo;\ X(t,\xi)=0\},\ \ t\ge0.$$
\end{theorem}

As mentioned earlier, Theorem \ref{t2.1} applies  to the
self-organized stochastic model \eqref{e1.2}, that is,
\begin{equation}\label{e2.8}
\barr{l} dX(t)-a\D\ {\rm
sign}\,(X(t)-X_c)dt=\sigma(X(t)-X_c)dW_t,\vsp X(0)=x-X_c\mbox{\ \
in }\calo.\earr
\end{equation}
If $x-X_c\ge0$ a.e. in $\calo$, then $X(t)-X_c\ge0$ a.e. on
$\calo$ for all $t\ge0$, $\pas$ and so, by Theorem \ref{t2.1},
it follows that $m(\calo\setminus\calo^t_0)\in L^1(0,\9)$, $\pas$,
which roughly speaking means that, "for almost all sequences"
$\{t_n\}\to\9$, we have $m(\calo\setminus\calo^t_0)\to0$, $\pas$

As regards the asymptotic result \eqref{e2.6}, one might expect
that $\ell=0$, $\pas$ Indeed, this is the case in the
deterministic case (see \cite{3}). For \eq\ \eqref{e2.1}, we have $$\lim_{t\to\9}X(t)=0\quad\mbox{in }L^1_{\rm loc}(\calo),$$ if the Gaussian noise $\sigma(X)W$ has
a finite number of modes,   that~is,
\begin{equation}\label{e2.9}
\sigma(X)W(t) =\sum^N_{k=1}\mu_ke_kX(t)\b_k(t)\ \mbox{on
}(0,\9)\times\calo\end{equation}
and
\begin{equation}\label{e2.10a}
\wt\mu(\xi)=\sum^N_{k=1}\mu^2_ke_k(\xi)>0,
\quad\ff\xi\in\calo.\end{equation}
More precisely,

\begin{theorem}\label{t2.2} Consider the situation of   Theorem {\rm\ref{t2.1}}. In addition, assume that \eqref{e2.9} holds and
$$\mu(t)=-\sum^N_{k=1}\mu_k e_k\b_k(t),\quad t\ge0.$$
Then:
\begin{itemize}
 \item[\rm(i)] $\dd\lim_{t\to\9}e^{\mu(t)}X(t)=0$ in $L^1(\calo)$, $\pas$ In particular, $\ell=0$ in the deterministic case $($cf. {\rm\cite{3}).}
     \item[\rm(ii)] If, additionally, \eqref{e2.10a} holds, then\vspace*{-5mm}\end{itemize}
\begin{equation}\label{e2.11}
\lim_{t\to\9} X(t)=0\quad\mbox{in }L^1_{\rm loc}(\calo).
\end{equation}

\quad Moreover,   for each compact set $K\subset\calo$,
\begin{equation}\label{e2.12}
\barr{r}\dd\int_K\!\!  X(t,\xi)d\xi\le
|x|_2 m(K)^{1/2}\exp\!\(\!\sup_{K}(\wt\mu)^{\!\!1/2}
\(\sum^N_{k=1}\b_k(t)\)^{\!\!\!1/2}\)e^{-\frac t2\,\inf\limits_{K'}\wt\mu}\!,\vsp t\ge0,\ \pas,
\earr\hspace*{-2mm}
\end{equation}
\begin{itemize}\item[\ ] where $K'\subset\calo$ is any compact neighborhood of $K$. In particular $($by the law of the iterated logarithm for Brownian motion$)$, there exists a constant $\rho_K>0$ such that, for $\mathbb{P}$-a.e. $\oo\in\ooo$\vspace*{-5mm}\end{itemize}
\begin{equation}\label{e2.12prim}
\int_KX(t,\xi,\oo)d\xi\le|x|_2m(K)^{1/2}e^{-\rho_Kt},\quad\ff t\ge t_0(\oo).\end{equation}
\end{theorem}

  We note that,
if $|\mu_1|>0$, then assumption \eqref{e2.10a} holds, because the
first eigenfunction of the Laplace operator is strictly positive
on $\calo$ (see, e.g.,\break \cite[p.~340]{Evans}).

The proofs of Theorems \ref{t2.1} and \ref{t2.2} are given in
Sections 3 and 4, respectively. For simplicity, we take $a=1$ in
\eqref{e2.1}.

\section{Proof of Theorem \ref{t2.1}}
\setcounter{equation}{0}

Consider the approximating \eq
\begin{equation}\label{e3.1}
\barr{l}
dX_\lbb(t)-\D(\psi_\lbb(X_\lbb(t))+\lbb X_\lbb(t))dt=\sigma(X_\lbb(t))dW_t\ \mbox{ in }(0,\9)\times\calo,\vsp
X_\lbb(0)=x\ \mbox{ on }\calo,\vsp
X_\lbb=0\ \mbox{ on }(0,\9)\times\pp\calo,\earr\end{equation}
where $\lbb\in(0,1)$ and
$$\psi_\lbb(r)=\dd\frac1r\ (1-(1+\lbb\psi)\1(r))=\left\{\barr{rl}
\dd\frac r\lbb\quad&\mbox{if }|r|\le\lbb,\vsp
1\quad&\mbox{if }r>\lbb,\vsp
-1\quad&\mbox{if }r<-\lbb.\earr\right.$$
(Here $1$ is the identity map.) As shown in \cite[Theorem 2.2]{4} and \cite[Proposition 3.5]{5}, \eq\ \eqref{e3.1} has a unique solution $X_\lbb$ in the sense of Definition \ref{d2.1} and
$$X_\lbb\in L^2(\ooo,C([0,T];H))\cap L^2(0,T;L^2(\ooo,H^1_0(\calo))).$$
Moreover, since $x\in L^4(\calo)$ and $x\ge0$, also $X_\lbb\ge0$ on $(0,\9)\times\calo\times\ooo$ and, as proved in \cite{5}, for $\lbb\to0$, we have for all $T>0$,
\begin{equation}\label{e3.2}
\barr{rcll}
X_\lbb&\to& X&\mbox{weakly in $L^2(\ooo\times(0,T)\times\calo),$}\vsp
&&&\mbox{weak$^*$ in $L^\9(0,T;L^2(\ooo,L^2(\calo))$ and}\vsp
&&&\mbox{strongly in $L^2(\ooo;C([0,T];H)),$}\earr\end{equation}
and
\begin{equation}\label{e3.3}
\barr{rcll}
\psi_\lbb(X_\lbb))+\lbb X_\lbb&\to&\eta&\mbox{weakly in }L^2(\ooo\times(0,T)\times\calo)\earr\end{equation}
\begin{equation}\label{e3.4}
\eta\in\psi(X)\mbox{\ \ a.e. on }\ooo\times(0,\9)\times\calo.\end{equation}
By Ito's formula and the monotonicity of $\psi_\lbb$, we have (cf. \cite[Proof of Theorem 2.8 and Remark 2.9(iii)]{8}
$$\barr{lcl}
\dd\frac12\ |X_\lbb(t)|^2_2&+&\dd\int^t_0\int_\calo\na(\psi_\lbb(X_\lbb)+\lbb X_\lbb)\cdot\na X_\lbb d\xi\,ds\vsp
&=&\dd\frac12\ |x|^2_2+\dd\frac12\sum^\9_{k=1}\mu^2_k\int^t_0\int_\calo|X_\lbb(s)e_k|^2d\xi\,ds\vsp
&+&
\dd\sum^\9_{k=1}\mu_k\int^t_0\int_\calo X^2_\lbb(s)e_kd\b_k(s),\  t\in[0,\9),\ \pas,\earr$$
and $t\to X_\lbb(t)\in L^2(\calo)$ is continuous $\pas$

This yields, since $\|e_k\|_\9\le C\lbb_k$, because $d\le 3$,
$$\barr{lcl}
|X_\lbb(t)|^2_2&+&2\lbb\dd\int^t_0|\na X_\lbb(s)|^2_2ds\le|x|^2_2+C\dd\sum^\9_{k=1}\mu^2_k\lbb^2_k\int^t_0|X_\lbb(s)|^2_2ds\vsp
&+&\dd\sum^\9_{k=1}\mu_k\int^t_0\int_\calo X^2_\lbb(s)e_kd\xi\,d\b_k(s).\earr$$
Then, by \eqref{e2.2} and by the Burkholder--Davis--Gundy inequality, we obtain the estimate that for some constant $C_T>0$ and all $\lbb\in(0,1)$
\begin{equation}\label{e3.5}
\mathbb{E}\sup_{t\in[0,T]}|X_\lbb(t)|^2_2+\lbb E\int^t_0|\na X_\lbb(s)|^2_2ds\le C_T|x|^2_2.\end{equation}
Now, arguing as in Proposition 3.5 in \cite{5}, we consider a \fc\ $\vf_\lbb\in C^3_b(\mathbb{R})$ such that $\vf_\lbb(0)=0$ and
\begin{equation}\label{e3.6}
\barr{l}
\vf'_\lbb(r)=\dd\frac r\lbb\ \mbox{for }|r|\le\lbb,\ \ \vf'_\lbb(r)=1+\lbb\ \mbox{for }r\ge2\lbb,\vsp
\vf'_\lbb(r)=-1-\lbb\ \mbox{for }r\le-2\lbb\mbox{ and }0\le\vf''_\lbb(r)\le\dd\frac C\lbb,\earr\end{equation}
for all $r\in \mathbb{R}$ and some $C>0$.

This is a smooth approximation of the \fc\ $r\to|r|$ and it is easily seen that
\begin{equation}\label{e3.6prim}
|\vf'_\lbb(r)-\psi_\lbb(r)|\le C\lbb,\ \ \ff r\in \mathbb{R},\ \lbb>0.\end{equation}
We set $Y^\vp_\lbb=(1+\vp A)\1 X_\lbb$ and note that
$$\barr{l}
dY^\vp_\lbb(t)+A(1+\vp A)\1(\psi_\lbb(X_\lbb(t))+\lbb X_\lbb(t))dt=(1+\vp A)\1\sigma(
X_\lbb(t))dW_t\vsp
Y^\vp_\lbb(0)=(1+\vp A)\1 X_\lbb(0),\ \ \vp>0.\earr$$
Also, the process $t\to Y^\vp_\lbb(t)$ is \ct\ $H^1_0(\calo)$-valued on $[0,T]$. Then, by Ito's formula applied to the $H^1_0$-valued process $Y^\vp_\lbb$, we have

\begin{equation}\label{e3.7}
\barr{l}
\dd\int_\calo\vf_\lbb(Y^\vp_\lbb(t,\xi))d\xi\vsp
\quad+
\dd\int^t_0\!\!\int_\calo\!\!\na((1+\vp A)\1(\psi_\lbb(X_\lbb(s,\xi){+}\lbb X_\lbb(s,\xi)){\cdot}\na\vf'_\lbb(Y^\vp_\lbb(s,\xi))d\xi\,ds\vsp
\quad=\dd\int_\calo\vf_\lbb((1+\vp A)\1x)d\xi\vsp
\quad+\dd\sum^\9_{k=1}\mu^2_k\dd\int^t_0\int_\calo\vf''_\lbb(Y^\vp_\lbb(s,\xi))|(1+\vp A)\1(X_\lbb e_k)(s,\xi)|^2ds\,d\xi\vsp
\quad+\dd\sum^\9_{k=1}\mu_k\dd\int^t_0\left<\vf'_\lbb(Y^\vp_\lbb(s,\xi)),
(1+\vp A)\1(X_\lbb e_k)(s,\xi))\right>_2d\b_k(s),\vsp \quad\hfill\ff t\ge0,\ \pas\earr\end{equation}
Now, recalling that $X_\lbb\in L^2(0,T;L^2(\ooo,H^1_0(\calo))$ for all $\lbb>0$, we have that $Y^\vp_\lbb\to X_\lbb$ strongly in $L^2(0,T;H^1_0(\calo)), \pas$ as $\vp\to0$. Similarly, for all $T>0$, we have

$$\barr{rcl}
\na((1+\vp A)\1(\psi_\lbb(X_\lbb)+\lbb X_\lbb))&\to&\na(\psi_\lbb(X_\lbb)+\lbb X_\lbb)\vsp
(1+\vp A)\1(X_\lbb e_k)&\to&X_\lbb e_k\earr$$
 strongly in $L^2(0,T;L^2\calo)),\ \pas$ as $\vp\to0$.

Furthermore, it is easy to see that by \eqref{e3.5} and the Burkholder--\break Davis--Gundy inequality for $p=1$, the stochastic term converges in\break $L^1(\ooo;C([0,T],L^2(\calo))$, as $\vp\to0$. Also, the first term in \eqref{e3.7} converges for a.e. $t\in[0,T]$ after passing to a subsequence $\vp_n\to0$. So, altogether, we obtain
\begin{equation}\label{e3.7a}
\barr{l}
\dd\int_\calo\vf_\lbb(X_\lbb(t))d\xi\vsp
\quad\dd+
\dd\int^t_0ds\int_\calo\na(\psi_\lbb(X_\lbb(s))+\lbb X_\lbb(s))\cdot\na\vf'_\lbb(X_\lbb(s))d\xi\vsp
\quad=\dd\int_\calo\vf_\lbb(x)d\xi+\sum^\9_{k=1}
\mu^2_k\int^t_0\int_\calo\vf''_\lbb(X_\lbb(s))
|X_\lbb(s)e_k|^2d\xi\,ds\vsp
\quad+\dd\sum^\9_{k=1}\mu_k\int^t_0\left<X_\lbb(s)e_k,
\vf'_\lbb(X_\lbb(s))\right>_2d\b_k(s)\mbox{ for a.e. }t>0,\ \pas\earr\end{equation}
On the other hand, by the $L^2(\calo)$-continuity of $X_\lbb$ it follows that the first term in \eqref{e3.7a} is \ct, as are all the other terms in \eqref{e3.7a}. Hence, \eqref{e3.7a} holds for all $t\ge0$, $\pas$

On the other hand, by \eqref{e3.6} we have the following estimate
\begin{equation}\label{e3.8}
\barr{l}
\dd\sum^\9_{k=1}\mu^2_k\int^t_0\int_\calo\vf''_\lbb(X_\lbb)|X_\lbb e_k|^2d\xi\,ds\vsp
\qquad\le4\lbb C\dd\sum^\9_{k=1}\mu^2_k\lbb^2_k\int^t_0\int_\calo\dubluunu_\lbb(s,\xi) d\xi\,ds,\ \pas,\earr\end{equation}
where $\dubluunu_\lbb$ is the characteristic \fc\ of the set $$\{(s,\xi,\oo)\in(0,\9)\times\calo\times\ooo;\ 0\le X_\lbb(s,\xi,\oo)\le2\lbb\}.$$ Now, we prove
\begin{equation}\label{e3.9}
\lim_{\lbb\to0}\int_\calo\vf_\lbb(X_\lbb(t,\xi))d\xi=\int_\calo X(t,\xi)d\xi,\ \ \ff t\ge0\mbox{ weakly in }L^2(\ooo).\end{equation}
Indeed, by \eqref{e3.6} we have for fixed $t\ge0$
$$\barr{lcl}
\dd\int_\calo\vf_\lbb(X_\lbb(t,\xi))d\xi&=&
(1+\lbb)\dd\int_{[X_\lbb(t,\xi)\ge2\lbb]}X_\lbb(t,\xi)d\xi\vsp
&&+\dd\frac1{2\lbb}\int_{[X_\lbb\le\lbb]}X^2_\lbb(t,\xi)d\xi+\dd\int_{[\lbb\le X_\lbb\le2\lbb]}\vf_\lbb(X_\lbb(t,\xi))d\xi.\earr$$
Taking into account \eqref{e3.2} and that $\vf_\lbb(r)\le C\lbb$ for $r\in[\lbb,2\lbb]$, this yields
\begin{equation}\label{e3.10}
\int_\calo\vf_\lbb(X_\lbb(t,\xi))d\xi=\int_\calo X_\lbb(t,\xi)d\xi+o(\lbb),\ \mbox{a.e. }\ff t\ge0,\ \pas\end{equation}
We also note that, by \eqref{e3.2}, we have
\begin{equation}\label{e3.11}
\barr{lcll}
X_\lbb(t)&\to&X(t)&\mbox{weakly in }L^2(\ooo\times\calo),\ \mbox{ for all }t>0.\earr\end{equation}
(Indeed, $\{X_\lbb(t)\}$ is strongly convergent to $X(t)$ in $L^2(\ooo;H)$ for each\break $t\in[0,\9)$ and is bounded in $L^2(\ooo\times\calo)$ for all $t\in[0,\9).$)

Then, by \eqref{e3.10} and \eqref{e3.11} we find that
\begin{equation*}\label{e3.12}
\lim_{\lbb\to0}\int_\calo\vf_\lbb(X_\lbb(t,\xi))d\xi=\int_\calo X(t,\xi)d\xi,\ \ \ff t\ge0\mbox{ weakly in $L^2(\ooo)$,}\end{equation*}as claimed.

Now, we set
\begin{eqnarray}
I_\lbb(t)&=&\int^t_0\int_\calo(\na\psi_\lbb(X_\lbb)+\lbb\na X_\lbb)\cdot\na\vf'_\lbb(X_\lbb)d\xi\,ds,\ t\ge0,\label{e3.13}\\[1mm]
M_\lbb(t)&=&\sum^\9_{k=1}\mu_k\int^t_0\left<X_\lbb e_k,\vf'_\lbb(X_\lbb)\right>_2d\b_k(s),\ t\ge0.\label{e3.14}
\end{eqnarray}

We recall  that, by (H2),
$$\barr{lcl}
M_\lbb(t)&=&\dd\int^t_0\<\vf'_\lbb(X_\lbb(s)),\sigma(X_\lbb(s))dW(s)\>_2\vsp
&=&\dd\int^t_0\<\sigma(X_\lbb(s))^*\vf'(X_\lbb(s)),dW(s)\>_2.\earr$$
where, for $h\in L^2(\calo)$,
$$\sigma(X_\lbb(s))h=\sum^\9_{k=1}\mu_k\<e_k,h\>_2X_\lbb(s)e_k.$$
We shall prove below that, for the adjoint operators $\sigma(X_\lbb(s))^*$ on $L^2(\calo)$ we have, for all $T>0$,
\begin{equation}\label{e3.14prim}
\sigma(X_\lbb)^*\vf'_\lbb(X_\lbb)\to\sigma(X)^*\eta\mbox{ weakly in }L^2((0,T)\times\calo\times\ooo)\mbox{ as }\lbb\to0.
\end{equation}
This implies that
\begin{equation}\label{e3.15a}
\barr{r}
\dd\lim_{\lbb\to0}M_\lbb(t)=M(t)=\sum^\9_{k=1}\mu_k
\int^t_0\<X(s)e_k,\eta\>_2d\b_k(s)\vsp\mbox{ weakly in }L^2(\ooo),\ \ff t\ge0.\earr\end{equation}
Now, let us prove \eqref{e3.14prim}. First, we note that by \eqref{e3.3}, \eqref{e3.4}, \eqref{e3.5} and \eqref{e3.6prim} as $\lbb\to0$
\begin{equation}\label{e3.15prim}
X_\lbb\to X\mbox{ and }\vf'_\lbb(X_\lbb)\to\eta\mbox{ weakly in }L^2((0,T)\times\calo\times\ooo),\end{equation}
and that, by \eqref{e3.6} and \eqref{e3.4},
\begin{equation}\label{e3.15secund}
|\vf'_\lbb(X_\lbb)|_\9,|\eta|_\9\le2,\end{equation}
where the norm refers to $L^\9((0,T)\times\calo\times\ooo).$ \eqref{e3.15secund} implies that, for some constant $C=C(T,\calo)>0$
$$\barr{l}
\mathbb{E}\dd\int^T_0|\sigma(X_\lbb(s))^*\vf'(X_\lbb(s))|^2_2ds\vsp
\quad\quad=\mathbb{E}\dd\int^T_0\sup_{|h|_2\le1}\<\vf'(X_\lbb(s)),
\sum^\9_{k=1}\mu_k\<e_k,h\>_2 X_\lbb(s)e_k\>^2_2ds\vsp
\quad\quad\le C\dd\sum^\9_{k=1}\mu^2_k\lbb^2_\lbb\mathbb{E}\dd\int^T_0|X_\lbb(s)|^2ds,\earr$$
which, by \eqref{e3.15prim} is uniformly bounded for $\lbb\in(0,1).$ Hence
\begin{equation}\label{e3.15tert}
\{\sigma(X_\lbb)^*\vf'(X_\lbb)\}_{\lbb\in(0,1]}\mbox{ is bounded in }L^2((0,T)\times\calo\times\ooo).\end{equation}
Now, let $F\in L^\9((0,T)\times\ooo;H^1_0(\calo)).$ Then
$$\barr{l}
\left|\mathbb{E}\dd\int^T_0\<F(s),
\sigma(X(s))^*\eta(s)-
\sigma(X_\lbb(s))^*\vf'_\lbb(X_\lbb(s))\>_2ds\right|\vsp
\quad\quad\le
\left|\mathbb{E}\dd\int^T_0\<F(s),\sigma(X(s))
(\eta(s)-\vf'_\lbb(X_\lbb(s)))\>_2ds\right|\vsp
\quad\quad+\left|\mathbb{E}\dd\int^T_0\<F(s),\sigma(X(s)-X_\lbb(s)) \vf'_\lbb(X_\lbb(s))\>_2ds\right|\vsp
\quad\quad\le\left|\dd\sum^N_{k=1}\mu_k
\mathbb{E}\dd\int^T_0
\<\<F(s), X(s)e_k\>_2  e_k,\eta(s)-\vf'_\lbb(X_\lbb(s))\>_2ds\right|\vsp
\quad\quad+ \mathbb{E}\dd\int^T_0
\left(\dd\sum^\9_{k=N+1}\mu^2_k
\<e_k,\eta(s)-\vf'_\lbb(X_\lbb(s))\>^2_2\right)^{\frac12}
|F(s)X(s)|_2ds  \vsp
\quad\quad+\left|\dd\sum^N_{k=1}\mu_k
\mathbb{E}\dd\int^T_0\<e_k,\vf'_\lbb(X_\lbb(s))\>_2
\<F(s)e_k,X(s)-X_\lbb(s)\>_2ds  \right|\vsp
\quad\quad+ \mathbb{E}\dd\int^T_0
\left(\dd\sum^\9_{k=N+1}\mu^2_k
\<e_k,\vf'_\lbb(X_\lbb(s))\>^2_2\right)^{\frac12}
|F(s)(X(s)-X_\lbb(s))|_2ds
\earr$$
of which the second and fourth term by \eqref{e3.15prim}, \eqref{e3.15secund} and the boundedness of $F$ converge  to zero uniformly in $\lbb\in(0,1)$ as $N\to\9$. By \eqref{e3.15prim}, the same is true for the first term for each fixed $N$ as $\lbb\to0$. Furthermore, the third term is up to a constant $C(T,\calo)>0$ dominated by
$$|F|_{L^\9((0,T)\times\calo;H^1_0)}\sum^N_{k=1}\mu_k\mathbb{E}\int^T_0|X(s)-X_\lbb(s)|^2_{-1}ds,$$
which, for each fixed $N$ as $\lbb\to0$, also converges to zero by \eqref{e3.2}. Hence, first letting $\lbb\to0$ and then $N\to\9$ and, using \eqref{e3.15tert}, we obtain \eqref{e3.14prim}.

Then, by \eqref{e3.7a}, \eqref{e3.8}, \eqref{e3.9} and \eqref{e3.15a}, we have
\begin{equation}\label{e3.16}
\int_\calo X(t,\xi)d\xi+\wt I(t)=\int_\calo x(\xi)d\xi+M(t),\ \ff t\ge0,\ \pas,\end{equation}
where
\begin{equation}\label{e3.17}
\wt I(t)=w-\lim_{\lbb\to0}I_\lbb(t),\ \ t\ge0,\end{equation}
and $w-\lim_{\lbb\to0}$ denotes   weak limit in $L^2(\calo)$.

We set$$Z(t)=\int_\calo X(t,\xi)d\xi,\ t\ge0.$$
We see that $Z$ is a nonnegative semimartingale with $EZ(t)<\9$, $\ff t\ge0$.

By \eqref{e3.5} and \eqref{e3.2} and  lower-semicontinuity, it follows that, for all $T>0$,
\begin{equation}\label{e3.18a}
\mathbb{E}\left[\sup_{t\in[0,T]}|X(t)|^2_2\right]<\9,\end{equation}
where we note that $\sup_{t\in[0,T]}|X(t)|^2_2={\rm ess\ sup}_{t\in[0,T]}|X(t)|^2_2$ since $\pas$\break $t\mapsto|X(t)|^2_2$ is lower-semi\ct\ by Definition \ref{d2.1}. The latter then together with \eqref{e3.18a} implies that  $\pas$ the \fc\ $t\to X(t)$ is weakly \ct\ in $L^2(\calo)$ on $[0,\9)$ and so the \fc\ $t\to Z(t)$ is $\pas$ \ct\ on $[0,\9)$.
Define
$$I(t):=Z(0)-Z(t)+M(t),\ \  t\ge0,$$
then $I$ is a \ct\ version of $\wt I$. We note that, clearly, by \eqref{e3.17} for all $0\le s\le t$
$$\wt I(s)\le\wt I(t),\ \pas$$
with the $\mathbb{P}$-exceptional set depending on $s,t$.

Hence (first considering all rational $s,t\in[0,\9)$, $0\le s\le t$), we conclude by continuity that
$$I(s)\le I(t),\ \ \ff\ 0\le s\le t,\ \pas,$$
i.e. $I$ is a  $\pas$ nondecreasing process.

Hence, altogether we have
$$Z(t)+I(t)=Z(0)+M(t),\ \ \ff t\ge0,$$
where $M$ is a \ct\ local martingale and $I$ is an a.s. nondecreasing process. Then, by \cite[p.~139]{7} we may conclude that
\begin{equation}\label{e3.18}
\exists\lim_{t\to\9}Z(t)<\9,\ I(\9)<\9,\ \pas\end{equation}
It follows therefore that there exists
\begin{equation}\label{e3.19}
\ell=\lim_{t\to\9}\int_\calo X(t,\xi)d\xi,\ \pas\end{equation}

Fix $t\ge0$. Noting that $\pas$
\begin{equation}\label{e3.27}
\barr{lcl}
 I_\lbb(t)&\ge&\dd\int^t_0\na\psi_\lbb(X_\lbb(s))
 \cdot\na\vf'_\lbb(X_\lbb(s))ds\vsp
 &=&\dd
\dd\int^t_0\na\psi_\lbb(X_\lbb(s))
\cdot\na\psi_\lbb(X_\lbb(s))ds,\earr\end{equation}
it follows by \eqref{e3.3}, \eqref{e3.5} and
\eqref{e3.17} that, as $\lbb\to0$,
\begin{equation}\label{e3.28}
\psi_\lbb(X_\lbb)\to\eta\mbox{ weakly in }L^2((0,T)\times\ooo;H^1_0(\calo)).\end{equation}
This, as well as \eqref{e3.17}, remains
true if $\mathbb{P}$ is replaced by $\rho\cdot\mathbb{P}$ for every $\rho\in L^\9(\ooo)$, $\rho\ge0$. Hence \eqref{e3.27} and \eqref{e3.17} imply
$$\mathbb{E}\left[\int^t_0|\na\eta|^2_2ds\ \rho\right]\le\liminf_{\lbb\to0}\mathbb{E}[I_\lbb(t)\rho]=\mathbb{E}[\wt I(t)\rho].$$
Since $\rho\in L^\9(\ooo)$, $\rho\ge0$, was arbitrary, this implies that
$$\int^t_0|\na\eta|^2_2ds\le\wt I(t),\ \pas$$
Hence, by continuity,
$$\int^t_0|\na\eta|^2_2ds\le I(t),\ \ \ff\ t\ge0,\ \pas$$
and, consequently, by \eqref{e3.18},
\begin{equation}\label{e3.29}
\lim_{t\to0}\int^t_0|\na\eta|^2_2ds\le I(\9)<\9,\ \pas\end{equation}
Now, by the Sobolev embedding theorem, we have by \eqref{e3.29}  that
\begin{equation}\label{e3.30}
\int^\9_0dt\(\int_\calo|\eta|^{p^*}d\xi\)^{\frac2{p^*}}<\9,\ \pas,\end{equation}
where $\frac1{p^*}=\frac12-\frac 1d$ for $d>2$, $p^*\in[2,\9)$ for $d=2$ and $p^*=\9$ for $d=1$.

Recalling that $\eta\in{\rm sign}\ X=1$ on $[X\ne0]$, a.e. on $(0,\9)\times\calo\times\ooo$, it follows by \eqref{e3.30} that
$$\int^\9_0(m(\calo\setminus\calo^t_0))^{\frac2{p^*}}dt<\9,\ \pas,$$
which implies \eqref{e2.7}, as claimed. This completes the proof of Theorem \ref{t2.1}.

\section{Proof of Theorem \ref{t2.2}}
\setcounter{equation}{0}

Assume in this section that \eqref{e2.9} holds. We recall that
\begin{equation}\label{e4.1}
\mu=-\sum^N_{k=1}\mu_k e_k\b_k,\ \ \wt\mu=\sum^N_{k=1}\mu^2_k e^2_k\end{equation}
and that the initial datum $x$ belongs to $L^4(\calo)$.

Take
\begin{equation}\label{e4.2}
Y(t)=e^{\mu(t)}X(t),\ \ \ff\ t\ge0.\end{equation} Then we
have (see \cite[Lemma~4.1]{5})
\begin{equation}\label{e4.3}
\barr{ll}
\dd\frac d{dt}\,Y(t)=e^\mu\D\psi(e^{-\mu}Y)-\dd\frac12\ \wt\mu Y,&\ff t\ge0,\ \pas,\vsp
Y(0)=x&\mbox{on }\calo,\vsp
\psi(e^{-\mu}Y)\in H^1_0(\calo),&\ff t\ge0,\ \pas,\earr\end{equation}
where the derivative $\frac d{dt}$ is taken in $H^{-1}(\calo)$.
(Recall that $\psi(r)={\rm sign}\ r$ and in \eqref{e4.3}, by Definition \ref{d2.1}, there arises a section $\eta$ of sign$\,(e^{-\mu}Y)$.)

First, we shall establish a  few estimates on the solution $Y$ to
\eqref{e4.3}, which have also an interest in themselves.

\begin{lemma}\label{l4.1} We have
\begin{equation}\label{e4.6}|Y(t)|_2\le|x|_2,\quad\ff t\ge0,\pas\end{equation}\end{lemma}

\noindent{\bf Proof.} Consider   the solution
$Y_\lbb$ to the approximating equation
\begin{equation}\label{e4.4}
\barr{l}
\dd\frac{dY_\lbb}{dt}=e^\mu\D(\psi_\lbb(e^{-\mu}Y_\lbb)+\lbb
e^{-\mu}Y_\lbb)-\frac12\ \wt\mu Y_\lbb,\vsp Y_\lbb(0)=X,\q
Y_\lbb\in L^2(0,T;H^1_0(\calo)),\earr\end{equation} which
corresponds to \eqref{e3.1}, i.e. $Y_\lbb=e^\mu X_\lbb$.  Multiplying \eqref{e4.4} by
 $Y_\lbb$  and integrating over $\calo$, we obtain
\begin{equation}\label{e4.5}\barr{r}
\dd\frac12\ \frac d{dt}\,|Y_\lbb(t)|^2_2+
\int_\calo\na(\psi_\lbb(e^{-\mu}Y_\lbb)+\lbb
e^{-\mu}Y_\lbb)\na(e^\mu Y_\lbb)d\xi\vsp\dd=-\frac12\int_\calo Y^2_\lbb\wt\mu d\xi   \le0,\mbox{ for a.e.
}t>0,\earr\end{equation} because $\wt\mu\ge0$, a.e. on
$\calo\times\ooo.$ On the other hand, recalling that
\begin{equation}\label{e4.4a}
\na(\psi_\lbb)(z)=\left\{\barr{cl} \dd\frac1\lbb\,\na z&\q\mbox{if
}|z|<\lbb,\vsp 0&\q\mbox{if }|z|\ge\lbb,\earr\right.\end{equation}
we get, by \eqref{e4.5},
$$\barr{lcl}
\dd\frac12\ \frac
d{dt}\,|Y_\lbb(t)|^2_2&+&\dd\frac{1}{\lbb}
\dd\int_\calo \dubluunu^{**}_\lbb
e^{2\mu}[|\na X_\lbb|^2+2X_\lbb\na X_\lbb\cdot\na \mu]d\xi\vsp
&+&\lbb\dd\int_\calo e^{2\mu}[|\na X_\lbb|^2+2X_\lbb\na X_\lbb\cdot\na \mu] d\xi\le0\quad\mbox{a.e. }t\ge0,\earr$$ where
 $\dubluunu^{**}_\lbb$  is the characteristic function
of $\{(t,\xi);\ 0\le(e^{-\mu}Y_\lbb)(t,\xi)\le\lbb\}.$ This yields
$$\barr{lcl}
\dd\frac d{dt}\,|Y_\lbb(t)|^2_2&\le&
2\dd\int_\calo\(\frac1\lbb\ \dubluunu^{**}_\lbb+\lbb\)|X_\lbb|^2 e^{2\mu}|\na\mu|^2 d\xi\vsp
&\le&2\lbb\dd\int_\calo(1+|X_\lbb|^2)e^{2\mu}|\na\mu|^2 d\xi,
  \q\mbox{a.e. }t>0.\earr$$

%%%%%%%%%%%

Integrating, we obtain
\begin{equation}\label{e4.6a}
|Y_\lbb(t)|^2_2\le|x|^2_2+2\lbb\dd\int^t_0|(1+X^2_\lbb)^{1/2}e^\mu|\na\mu||^2_2ds,\q\ff t\ge0,\ \pas
\end{equation}
Defining $Y^{(N)}_\lbb:=X_\lbb(e^\mu\wedge N)$ and $Y^{(N)}:=X(e^\mu\wedge N)$, $N\in\mathbb{N}$, we deduce from \eqref{e3.2} that for all $\rho\in L^\9(\ooo)$, $\rho\ge0$,  as  $\lbb\to\9,$
$$Y^{(N)}_\lbb\to Y^{(N)}\mbox{ weak}^*\mbox{ in }L^\9(0,T;L^2(\ooo,\rho\mathbb{P};L^2(\calo)).$$
Hence
\begin{equation}\label{e4.6secund}
\st{t\in[0,T]}{\rm ess\ sup}  \mathbb{E}[|Y^{(N)}(t)|^2_2\rho]\le\liminf_{\lbb\to\9}\ \st{t\in[0,T]}{\rm ess\ sup}\mathbb{E}[|Y^{(N)}_\lbb(t)|^2_2\rho]. \end{equation}
But, by \eqref{e4.6a}, for all $N\in\mathbb{N}$,
$$\mathbb{E}[|Y^{(N)}_\lbb(t)|^2_2\rho
]\le|x|^2_2
\mathbb{E}[\rho]+2\lbb\|\rho\|_\9C,\q\ff t\ge0,\ \pas,$$
where
$$C:=\int^T_0(\mathbb{E}|e^\mu|\na\mu||^4_4)^{1/2}dt\cdot\sup_{\lbb\in(0,1)}\
\st{t\in[0,T]}{\rm ess\ sup}(\mathbb{E}|1+X_\lbb|^4_4)^{1/2}
$$
is finite by \cite[Lemma 3.1]{5}. Hence, letting first $\lbb\to0$ and then $N\to\9$ in \eqref{e4.6secund}, since $\rho\in L^\9(\ooo)$, $\rho\ge0$, was arbitrary, we obtain that
$$|Y(t)|^2_2\le|x|^2_2\q\mbox{for a.e. }t>0,\ \pas$$
Now \eqref{e4.6} follows, since $\pas$ $t\to|Y(t)|^2_2$ is lower-semi\ct.\bigskip

Now, let us turn to the proof of Theorem \ref{t2.2}.

To prove (i), let us assume that for some sequence $t_n\to\9$ we have that
\begin{equation}\label{e4.7}
|Y(t_n)|_1\ge\de>0,\q\ff n\in\mathbb{N}.\end{equation}
Here and below $Y(t)=Y(t,\oo)$ for a fixed $\oo\in\ooo$ such that \eqref{e4.6} holds. By \eqref{e4.6}, selecting a subsequence if necessary, we have $Y(t_n)\to g$ weakly in $L^2(\calo)$ as $n\to\9$. We have that $g\ge0$ and by \eqref{e4.7}
\begin{equation}\label{e4.8}
g\not\equiv0.\end{equation}
We recall from the proof of Theorem \ref{t2.1} that $t\mapsto\int_\calo X(t)d\xi$ is \ct, hence so is $t\to\int_\calo Y(t)d\xi.$ So, for every $n\in\mathbb{N}$, there exists $\vp_n>0$ such that
\begin{equation}\label{e4.9}
\left|\int_\calo Y(t)d\xi-\int_\calo Y(t_n)d\xi\right|\le\frac1n,\q\ff t\in(t_n-\vp_n,t_n+\vp_n).\end{equation}
It follows by \eqref{e2.7} that for some subsequence $t_{n_k}\to\9$ there exist $s_k\in(t_{n_k}-\vp_{n_k},t_{n_k}+\vp_{n_k})$, $k\in\mathbb{N},$ such that
$$\int\dubluunu_{\{X(s_k)\ne0\}}d\xi=
m(\calo\setminus\calo^{s_k}_0)\to0\mbox{\ \ as }k\to\9.$$
Hence, selecting another subsequence if necessary, we have
$$\dubluunu_{\{X(s_k)\ne0\}}\to0\mbox{\ \ a.e. as }k\to\9$$
and by \eqref{e4.6} that  $X(s_k)\to\wt g$ weakly in $L^2(\calo).$

As a consequence of the first, we obtain
$$Y(s_k)=Y(s_k)\dubluunu_{\{X(s_k)\ne0\}}\to0\mbox{\ \ a.e. as }k\to\9,$$
which, in turn, implies that $\wt g=0$. Hence, by \eqref{e4.9}
$$\int g\,d\xi=\lim_{k\to\9}\int Y(t_{n_k})d\xi=\lim_{k\to\9}\int Y(s_k)d\xi=\int\wt g\,d\xi=0.$$
Hence, $g=0$ a.e., since $g\ge0$. This contradiction to \eqref{e4.8} proves that a sequence $t_n\to0$ with \eqref{e4.7} does not exist and assertion (i) follows.

%%%%%%%%%%%%%

Clearly, to prove (ii), it suffices to prove  the exponential  decay
part of Theo\-rem \ref{t2.2} (ii). So, additionally, assume that \eqref{e2.10a} holds and let $K\subset\calo$, $K$ compact, and $K'\subset\calo$ a compact neighborhood of $K$, i.e., $K\subset \bld\circ{K'}$. Let $\mu^*\in C^\9_0(\calo)$ such that $0\le\mu^*\le1$, $\mu^*=1$ on $K$ and $\mu^*=0$ on $\calo\setminus K'$. Furthermore, let $C_K:=\inf_{K'}\wt\mu$. We multiply equation
\eqref{e4.4} by  $\mu^* Y_\lbb$   and integrate
over~$\calo$ to obtain

\begin{equation}\label{e4.12a}
\barr{l} \dd\frac12\ \frac
d{dt}\,|(\mu^*)^{\frac12}Y_\lbb(t)|^2_2+\dd\frac{C_K}2\,|(\mu^*)^{\frac12}Y_\lbb(t)|^2_2\vsp
\qquad\le\dd\frac12\ \frac d{dt}\ |(\mu^*)^{\frac12}Y_\lbb(t)|^2_2
+\dd\frac12\int_\calo\wt\mu\mu^* Y^2_\lbb d\xi\vsp
\qquad=-\dd\int_\calo\na(\psi_\lbb(e^{-\mu}Y_\lbb))\cdot \na(e^\mu \mu^*
Y_\lbb)d\xi\vsp
\qquad-\lbb\dd\int_\calo\na(e^{-\mu}Y_\lbb)\cdot\na(e^\mu  \mu^*
Y_\lbb)d\xi\vsp
\qquad=-\dd\frac1\lbb\int_\calo\dubluunu^{**}_\lbb\left[|\na X_\lbb|^2+2X_\lbb\na X_\lbb\cdot\(\na\mu+\frac12\ \frac{\na\mu^*}{\mu^*}\)\right]e^{2\mu}\mu^* d\xi\vsp
\qquad-\lbb\dd\int_\calo\left[|\na X_\lbb|^2+2X_\lbb\na X_\lbb\cdot\(\na\mu+\frac12\ \frac{\na\mu^*}{\mu^*}\)\right] e^{2\mu}\mu^*d\xi\vsp
\qquad\le\dd\int_\calo\(\frac1\lbb\ \dubluunu^*_\lbb+\lbb\)X^2_\lbb\left[2|\na\mu|^2\mu^*+\frac12\ \frac{|\na\mu^*|^2}{\mu^*}\right] e^{2\mu}d\xi\vsp
\qquad\le\dd2\lbb\int_\calo(1+X^2_\lbb)[|\na\mu|^2\mu^*+|\na(\mu^*)^{\frac12}|^2]e^{2\mu}d\xi
.\earr\end{equation}
Denoting the latter by $\frac\lbb2\ \eta_\lbb(t)$, we deduce that
$$\frac d{dt}\,(|(\mu^*)^{\frac12}Y_\lbb(t)|^2_2e^{C_Kt})\le\eta_\lbb(t)e^{C_Kt},\mbox{ for a.e. }t>0,\ \pas$$
Integrating from $0$ to $t$, we obtain
\begin {equation}\label{e4.14prim}
|(\mu^*)^{\frac12}Y_\lbb(t)|^2_2\le|(\mu^*)^{\frac12}x|^2_2 e^{-C_Kt}+\lbb\int^t_0 e^{C_K(s-t)}\eta_\lbb(s)ds,\ t\ge0,\ \pas\end{equation}
Now analogous arguments as in the proof of Lemma \ref{l4.1} imply that after letting $\lbb\to0$, inequality \eqref{e4.14prim} turns into
\begin{equation}\label{e4.14secund}
\barr{lcl}
|(\mu^*)^{\frac12}Y(t)|^2_2&\le&e^{-C_Kt}|(\mu^*)^{\frac12}x|^2_2,\quad t\ge0,\ \pas,\vsp
&\le&e^{-C_Kt}|x|^2_2.\earr\end{equation}
Hence
$$\barr{lcl}
\dd\int_KX(t)d\xi&=&\dd\int_KY(t)e^{\mu(t)}d\xi\vsp
&\le&|(\mu^*)^{\frac12}Y(t)|_2\(\dd\int_K
\exp\(2(\wt\mu)^{\frac12}\(\sum^N_{k=1}\b_k(t)^2\)^{\frac12}\)d\xi\)^{\frac12}\vsp
&\le&e^{-\frac{C_K}2\,t}|x|_2\exp\(\dd\sup_K (\wt\mu)^{\frac12}\(\sum^N_{k=1}\b_k(t)^2\)^{\frac12}\)
m(K)^{\frac12},\earr$$
i.e. \eqref{e2.12} is proved.

%%%%%%%%%%%%%%

\begin{remark} \label{r4.2} {\rm For existence of solutions to equation \eqref{e2.1} in the special case  \eqref{e2.9}, it is not absolutely necessary to assume that $\{e_k\}\subset H^1_0(\calo)$ is a basis of eigenfunctions for $A$. It suffices to assume that $e_k\in C^2(\ov\calo)$ and the proof of Theorem \ref{t2.2} is essentially the same. Then one might choose $e_k$, $1\le k\le N$, such that
$$\inf\{\wt\mu(\xi);\ \xi\in\ov\calo\}=\rho>0$$
and, in this case, the exponential decay   in Theorem \ref{t2.2}
is global in $\calo$. More precisely, in \eqref{e2.11} the compact sets $K$ and $K'$ can be replaced by $\calo$, and, in this case, \eqref{e2.12} strengthens to
\begin{equation}\label{e2.10prim}
\lim_{t\to\9}X(t)=0\quad\mbox{in }L^1(\calo),\ \pas,\end{equation}
and, therefore, $\ell=0$, $\pas$  The details are omitted.}\end{remark}

\begin{remark} \label{r4.3} {\rm If condition \eqref{e2.10a} does not hold, the following slightly weaker statements still hold. Since $\wt\mu$ is analytic on
$\calo$, the set $\{\xi_j\in\calo;\ \wt\mu(\xi_j)\}$ is countable and,
therefore, $\wt\mu(\xi)\ge\rho_K>0$, $\ff\xi\in K$,
for any compact $K\subset\calo\setminus\{\xi_j\}$. Then the proof
of Theorem \ref{t2.2} applies word by word and we have \eqref{e2.12} and \eqref{e2.12prim} in this case,
too.}\end{remark}

\end{document}